\documentclass[11pt]{article}
\usepackage[dvips]{graphics}
\input{amssym.def}
\input{amssym.tex}

\title{\Large\bf 
Continuum mechanics model of graphene\\ as a doubly-periodic perforated thin elastic plate}
\author{Yuri A.  Antipov\\
yantipov@lsu.edu\\ 
Department of Mathematics, Louisiana State University\\
Baton Rouge LA 70803, USA}

\date{}

\newcommand{\I}{\mathop{\rm Im}\nolimits}
\newcommand{\R}{\mathop{\rm Re}\nolimits}

\newcommand{\beqa}{\begin{eqnarray}}
\newcommand{\eeqa}[1]{\label{#1}\end{eqnarray}}
\newcommand{\bequ}{\begin{equation}}
\newcommand{\eequ}[1]{\label{#1}\end{equation}}

\newcommand{\Md}{\partial}

\newcommand{\ov}[1]{\overline{#1}}

\newcommand{\Ga}{\alpha}
\newcommand{\Gb}{\beta}
\newcommand{\Gd}{\delta}

\newcommand{\Gvf}{\varphi}
\newcommand{\Gg}{\gamma}
\newcommand{\Gc}{\chi}

\newcommand{\Gk}{\kappa}

\newcommand{\Gl}{\lambda}

\newcommand{\Gt}{\theta}

\newcommand{\Gr}{\rho}

\newcommand{\Gs}{\sigma}

\newcommand{\Go}{\omega}

\newcommand{\Gz}{\zeta}
\newcommand{\GD}{\Delta}
\newcommand{\GF}{\Phi}
\newcommand{\GG}{\Gamma}
\newcommand{\GL}{\Lambda}

\newcommand{\GY}{\Psi}

\newcommand{\CA}{{\cal A}}

\newcommand{\CD}{{\cal D}}

\newcommand{\CN}{{\cal N}}

\def\half{{\scriptstyle{1\over 2}}}

\newcommand{\beq}{\begin{equation}}
\newcommand{\eeq}{\end{equation}}
\newcommand{\barr}{\begin{eqnarray}}
\newcommand{\earr}{\end{eqnarray}}
\newcommand{\beqn}{\begin{equation*}}
\newcommand{\eeqn}{\end{equation*}}
\newcommand{\barrn}{\begin{eqnarray*}}
\newcommand{\earrn}{\end{eqnarray*}}

\newcommand{\fr}{\frac}

\begin{document}
\maketitle

\begin{abstract}

In this paper, a continuum mechanics model of graphene is proposed, and its analytical solution is derived.
Graphene is modeled as a doubly-periodic  thin elastic plate with a hexagonal cell having a circular hole at the hexagon center. 
Graphene is characterized by a general chiral vector and is subject to remote tension.
For the solution,  the Filshtinskii solution obtained for the symmetric case  is generalized for
any chirality. The method uses the doubly-periodic Kolosov-Muskhelishvili complex potentials, the theory
of the  elliptic Weierstrass function and quasi-doubly-periodic meromorphic functions and reduces the model
to an infinite system of linear algebraic equations with complex coefficients. Analytical expressions and numerical
values for the stresses are displacements are obtained and discussed. The displacements expressions possess the Young modulus and
Poisson ratio of the graphene bonds. They are derived  as functions of the effective graphene moduli available in the literature.

\end{abstract}


\setcounter{equation}{0}

\section{Introduction} 

Characterization of the mechanical properties of graphene and other  two-dimensional  materials is a challenging problem for their thickness. 
Different experimental and numerical simulation techniques are used to determine  elastic properties and in particular, estimate 
the Young modulus $E$ and Poisson ratio $\nu$ of graphene. By applying a force $F$ and measuring the relative elongation $\GD L$ for an equilibrium length $L_0$, the cross-section $A$, and the relative change in radius $\GD R/R$ the constants $E$ and $\nu$ 
of a single-walled carbon nanotube
were calculated \cite{lier} as $\fr{F}{A}/\fr{\GD L}{L_0}$ and $-\fr{\GD R}{R}/\fr{\GD L}{L_0}$, respectively.   In a similar manner these constants were recovered \cite{lier}
for graphene.
A thorough comparative analysis of different experimental measurements, {\it ab initio} calculations and computations based on density functional theory presented in  \cite{leb}, \cite{mem} shows that further work is needed to improve the accuracy of the data.
The finite element method was applied  \cite{cao} to investigate the accuracy of Raman spectroscopy method for measurement of the graphene elastic moduli. It was found that 
the determination of the pre-strain of a freestanding graphene membrane is essential for improvement of estimates of the elastic moduli 
determined via Raman spectroscopy. Continuum elasticity and tight-binding atomistic simulations were combined \cite{cad} to derive constitutive nonlinear stress-strain relations for graphene. These relations were employed to estimate the elastic moduli of graphene.
The majority of measurements and numerical calculations report that the value of graphene Young's modulus is $E\sim 1.0$ TPa \cite{liu}, \cite{mem}.
For the graphene Poisson ratio the data deviate significantly, $0.12\le \nu\le 0.413$ \cite{mem}. 

The values reported in the literature are effective elastic constants of graphene and  do not distinguish the bond moduli from
those obtained for a sample graphene membrane. To have estimates for the elastic moduli of the graphene bonds is important for mechanical models of graphene
when it is of interest to evaluate the displacements and strains and also  to analyze the mechanics of the defective van der Waals 
structures such as  hexagon–heptagon defects in graphene \cite{zha}.

The main aims of this paper  are (1) to develop a continuum mechanics model of graphene as a thin elastic doubly-periodic perforated hexagonal
plate and (2) to estimate the Young modulus and Poisson ratio of the graphene bonds as functions of the corresponding graphene effective  
moduli available in the literature.
In Section 2  we formulate the mathematical model of plane strain elasticity for a doubly-periodic plane with a hexagonal cell having
a circular hole in the center of the cell. The solution to the model problem is represented as a sum of the solution to the classical Kirsch problem [8] and the solution to an auxiliary problem with vanishing stresses at infinity and inhomogeneous boundary conditions 
on the holes boundary.  In Section 3 we employ the  Filshtinskii method \cite{fil} of elliptic functions to reduce the problem 
to an infinite system. The methodology is based on the use of the Kolosov-Muskhelishvili potentials [9] and their series representations  
in terms of the  Weierstrass elliptic function [10] and its derivatives and the Natanzon quasi-doubly-periodic function [11] and its derivatives.
At the final stage the method requires to satisfy the boundary conditions on the hole boundary which lead to an infinite system 
of linear algebraic equations for the series coefficients.
In Ref \cite{fil}  the infinite system was solved for the symmetric case (in the modern terminology, for the armchair - zigzag
symmetry) when the complex Kolosov-Muskhelishvili potentials satisfy the symmetry conditions
\beq
\ov{\GF(z)}=\GF(\bar z), \quad \ov{\GY(z)}=\GY(\bar z).
\label{1.1}
\eeq
In this case a unit cell is an equilateral triangle, the system is simplified, and the expansion coefficients are real.
We generalize the method to structures characterized by any chiral vector. In the general case  the condition (\ref{1.1}) is not valid, 
the unit cell is a hexagon, and  the coefficients are complex and  solve an infinite system with complex entries.
In Section 4 we derive series representations of the stresses and displacements in the perforated plane. The displacements
formulas possess the shear modulus $G=\fr{E}{2(1+\nu)}$ and the Poisson ratio $\nu$ of the graphene bonds (the material of the perforated plane) which are not available in the 
literature.
We modify the method \cite{fil} used for the derivation of the effective moduli $E^\circ$ and $\nu^\circ$  in terms of  $E$ and $\nu$
in the case of a triangular equilateral cell by reversing it and  derive  the elastic moduli  $E$ and $\nu$ of the graphene  bonds  as functions of the effective
moduli $E^\circ$ and $\nu^\circ$ measured by tests.
Finally, we report and discuss  our numerical results for the stresses, displacements and elastic moduli.

\setcounter{equation}{0}

\section{Formulation}\label{form}

\begin{figure}[t]
\centerline{
\scalebox{0.4}{\includegraphics{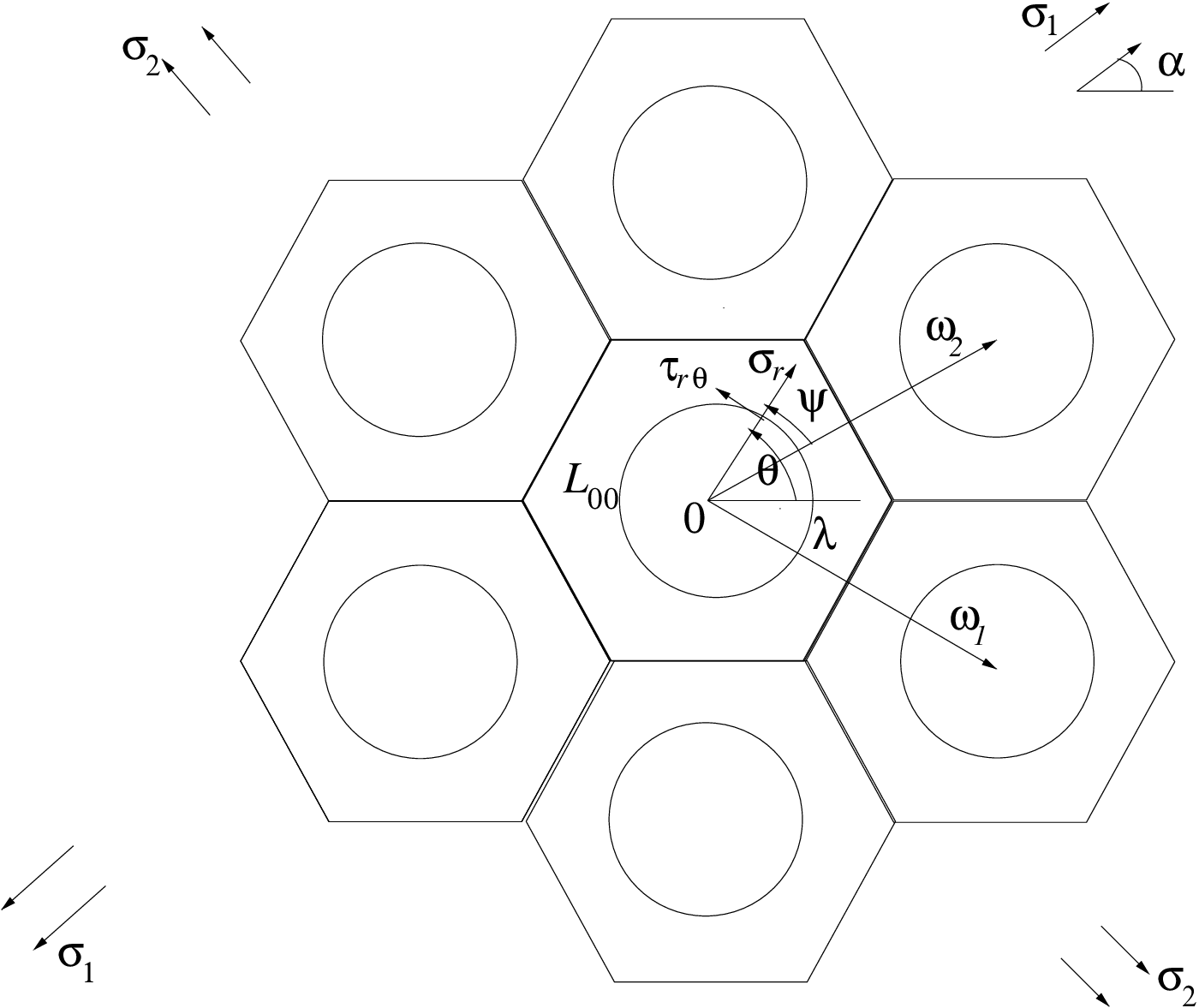}}}
\caption{Geometry of the problem
}
\label{fig1}
\end{figure}

Consider extension of an infinite single-layer graphene sheet by remote constant loads $\Gs_1$ and $\Gs_2$ (Fig. 1).
The lattice vectors of graphene are selected as 
\beq
\Go_1=\left(\fr{a\sqrt{3}}{2},-\fr{a}{2}\right), \quad \Go_2=\left(\fr{a\sqrt{3}}{2},\fr{a}{2}\right),
\label{2.1}
\eeq
where $a=|\Go_1|=|\Go_2|=142 \sqrt{3}=246$ pm is the lattice graphene constant and $\fr{a}{\sqrt{3}}=142$ pm is the nearest-neighbor distance between two carbon atoms.
The chiral vector $C_h$ \cite{sai} is expressed through the lattice vectors as 
\beq
C_h=m\Go_1+n\Go_2\equiv (m,n), \quad m,n=\pm 1,\pm 2,\ldots.
\label{2.1'}
\eeq
The armchair and zigzag symmetries correspond to the chiral vectors $(n,n)$ and $(n,0)$, respectively.  
 The graphene sheet
is modeled as a hexagonal  doubly-periodic perforated thin  elastic  plate. The loads $\Gs_1$ and $\Gs_2$  in the upper half-plane
make angles $\Ga\in [0,\pi/2]$ and $\Ga+\pi/2$ with the positive $x$-axis, respectively, and
the angle $\Ga$ in the range $0\le\Ga\le\fr{\pi}{6}$ relates to the chiral vector as
\beq
\Ga=\fr{\pi}{6}-\cos^{-1}\fr{C_h\cdot\Go_2}{|C_h| |\Go_2|}=\fr{\pi}{6}-\cos^{-1}\fr{2n+m}{2\sqrt{n^2+m^2+nm}}.
\label{2.1''} 
\eeq
The radius of the holes is denotes by $\Gl$, and
the boundary of the circular hole in each hexagonal cell or, equivalently, the graphene bonds, are free of traction.

The solution to the model problem, call it $P^t$, is comprised of the solution of the classical Kirsch problem  \cite{kir},  $P^K$,  of extension of an elastic plane with a 
circular hole by remote loads $\Gs_1$ and $\Gs_2$  and the plane problem $P$ of a doubly-periodic perforated plane (Fig. 1) with stresses vanishing at infinity
and  inhomogeneous boundary conditions on the boundary of the holes. These conditions are   chosen such that the total traction in Problem $P^t$ vanishes on the boundary.
The solution to Problem $P^K$ written in terms of the stresses $\Gs_r$ and $\tau_{r\Gt}$ in polar coordinates $(r,\Gt)$  has the form 
\beq
\Gs_r^K=\Gs_++\Gs_-\cos 2\psi, \quad \tau_{r\Gt}^K=-\Gs_-\sin 2\psi,
\label{2.2}
\eeq
where
\beq
\Gs_+=\fr{\Gs_1+\Gs_2}{2}, \quad \Gs_-=\fr{\Gs_1-\Gs_2}{2}, \quad \psi=\Gt-\Ga.
\label{2.3}
\eeq

Problem $P$ is governed by the equilibrium equations of plane elasticity 
\beq
\fr{\Md \Gs_x}{\Md x}+\fr{\Md \tau_{xy}}{\Md y}=0,
\quad \fr{\Md \tau_{xy}}{\Md x}+\fr{\Md \Gs_y}{\Md y}=0, \quad (x,y)\in \CD,
\label{2.4}
\eeq
and the compatibility conditions
\beq
\GD(\Gs_x+\Gs_y)=0,  \quad (x,y)\in \CD,
\label{2.5}
\eeq
where $\CD$ is the graphene domain (the perforated plane) and $\Gs_x, \Gs_y$, and $\tau_{xy}$ are the stress tensor components in the Cartesian coordinates $x,y$.
At infinity, the stresses vanish, while on the circle  $L_{mn}$, the boundary of the hole $D_{mn}$, the traction components have the values
\beq
\Gs_r=-\Gs_+-\Gs_-\cos 2\psi,\quad \tau_{r\Gt}=\Gs_-\sin 2\psi, \quad r=\Gl, \quad 0\le \Gt\le 2\pi.
\label{2.6}
\eeq
In addition, the stresses  $\Gs_x, \Gs_y$, and $\tau_{xy}$ have to be doubly-periodic functions and satisfy the condition
\beq
(\Gs_x,\Gs_y,\tau_{xy})(z)=(\Gs_x,\Gs_y,\tau_{xy})(z+m\Go_1 +n\Go_2), \quad m,n=0,\pm 1, \pm 2,\ldots,
\label{2.7}
\eeq
where 
\beq
\Go_1=\fr{a\sqrt{3}}{2}-\fr{ia}{2}, \quad \Go_2=\fr{a\sqrt{3}}{2}+\fr{ia}{2}
\label{2.8}
\eeq
are the lattice vectors written in the complex form and $z=x+iy\in H_{00}\setminus D_{00}$ is  a point in the unit hexagonal cell $H_{00}$ centered at the origin of the hole $D_{00}$
with the boundary $L_{00}$.

\setcounter{equation}{0}
  
\section{Solution of the boundary-value problem}\label{bvp} 

\subsection{Series representation of the solution in terms of the Weierstrass  and Natanzon functions}

The first step of the solution procedure is  to represent the stresses  by means  of the classical Kolosov-Muskhelishvili potentials 
$\GF(z)$ and $\Psi(z)$ analytic everywhere in the domain $\CD$ \cite{mus}
$$
\Gs_x+\Gs_y=2[\GF(z)+\ov{\GF(z)}]=4\R \GF(z),
$$
\beq
\Gs_y-\Gs_x+2i\tau_{xy}=2[\bar z\GF'(z)+\Psi(z)], \quad z\in \CD.
\label{3.1}
\eeq
As it is known, the stresses expressed through analytic functions by these relations automatically
satisfy the equilibrium equations (\ref{2.4}) and the compatibility condition (\ref{2.5}). Because of the structure and the applied  load symmetries these functions
have to be even,
\beq
\GF(z)=\GF(-z), \quad \Psi(z)=\Psi(-z), \quad z\in\CD.
\label{3.2}
\eeq
Also, due to the doubly-periodicity property (\ref{2.7}) of  the stresses and the representations (\ref{3.1}),  the function  $\GF(z)$ has to be doubly-periodic,
while the function $\Psi(z)$ is to be quasi-doubly-periodic,
\beq
\GF(z+\Go_j)-\GF(z)=0, \quad \Psi(z+\Go_j)-\Psi(z)=-\ov{\Go_j}\GF'(z), \quad j=1,2.
\label{3.3}
\eeq
General series representations of such functions has the form \cite{fil}
$$
\GF(z)=\Ga_0+\sum_{k=1}^\infty\Ga_{k}\fr{\Gl^{2k}\wp^{(2k-2)}(z)}{(2k-1)!},
$$
\beq
\Psi(z)=\Gb_0+\sum_{k=1}^\infty\Gb_{k}\fr{\Gl^{2k}\wp^{(2k-2)}(z)}{(2k-1)!}
-\sum_{k=1}^\infty\Ga_{k}\fr{\Gl^{2k}\CN^{(2k-1)}(z)}{(2k-1)!}.
\label{3.4}
\eeq
Here, $\Ga_k$ and $\Gb_k$ are arbitrary complex coefficients to be determined,  $\wp(z)$ is  Weierstrass' elliptic function \cite{han}
\beq
\wp(z)=\fr{1}{z^2}+\sum_{m,n}^{} {}'\left[\fr{1}{(z-w)^2}-\fr{1}{w^2}\right]
\label{3.5}
\eeq
and $\CN(z)$ is Natanzon's meromorphic function \cite{nat}
\beq
\CN(z)=\sum_{m,n}^{} {}'\bar w\left[\fr{1}{(z-w)^2}-\fr{2z}{w^3}-\fr{1}{w^2}\right].
\label{3.6}
\eeq
The accent on the summation sign indicates that the value $w=0$ is excluded from the summation, while $w$ is defined by
\beq
w=m\Go_1+n\Go_2, \quad m,n=0,\pm 1,\pm 2,\ldots.
\label{3.7}
\eeq
The series expansions  (\ref{3.4}) of the functions $\GF(z)$ and $\Psi(z)$ in terms of the even and odd order derivatives of the functions $\wp(z)$
and $\CN(z)$ converge uniformly and absolutely in the domain $\CD$. The derivatives  $\wp^{(2k-2)}(z)$ and $\CN^{(2k-1)}(z)$ are holomorphic
in the domain $\CD$ (all their poles lie in the holes that is in the exterior of the domain $\CD$).
The Weierstrass function $\wp(z)$ and all its derivatives 
\beq
\wp^{(k)}(z)=(-1)^k (k+1)!\left[\fr{1}{z^{k+2}}+\sum_{m,n}^{} {}'\fr{1}{(z-w)^{k+2}}\right]
\label{3.7'}
\eeq
are doubly-periodic,
\beq
\wp^{(k)}(z+\Go_j)-\wp^{(k)}(z)=0, \quad  j=1,2, \quad k=0,1,\ldots,
\label{3.8}
\eeq
while the Natanzon function and its derivatives 
$$
\CN'(z)=-2\sum_{m,n}^{} {}' 
\left[\fr{\bar w}{(z-w)^3}+\fr{\bar w}{w^3}\right], 
$$
\beq
\CN^{(k)}(z)=(-1)^k (k+1)!\sum_{m,n}^{} {}' \fr{\bar w}{(z-w)^{k+2}}, \quad k=2,3,\ldots,
\label{3.8'}
\eeq
are quasi-doubly-periodic,
\beq
\CN^{(k)}(z+\Go_j)-\CN^{(k)}(z)=\bar \Go_j\wp^{(k)}(z)+\Gg_j\Gd_{k0},  \quad  j=1,2, \quad k=0,1,\ldots,
\label{3.9}
\eeq
where $\Gd_{k0}$ is the Kronecker symbol, the constants $\Gg_1$ and $\Gg_2$  satisfy the following relations \cite{fil}:
\beq
\Gg_2\Go_1-\Gg_1\Go_2=\Gd_1\bar\Go_2-\Gd_2\bar\Go_1,  \quad \Gd_j=2\Gz\left(\fr{\Go_j}{2}\right), \quad j=1,2,
\label{3.9'}
\eeq
$\Gz(z)$ is the Weierstrass $\Gz$-function, $\wp(z)=-\Gz'(z)$, 
\beq
\Gz(z)= \fr{1}{z}+\sum_{m,n}^{} {}'\left[\fr{1}{z-w}+\fr{1}{w}+\fr{z}{w^2}\right],
\label{3.10}
\eeq
and $\Gd_1$ and $\Gd_2$ are cyclic constants. The $\Gz$-function is quasi-doubly-periodic with the property  \cite{han}
\beq
\Gz(z+\Go_j)-\Gz(z)=\Gd_j,  \quad j=1,2.
\label{3.11}
\eeq

We remark that  the derivatives of even order of the Weierstrass function and odd order of the Natanzon function are even.
Therefore the functions in (\ref{3.4}) are even.

\subsection{Determinafiltion of the coefficients $\Ga_0$ and $\Gb_0$}

We proceed now with determination of the coefficients $\Ga_k$ and $\Gb_k$ ($k=0,1,\ldots.$)
Consider any two congruent points $A\in\CD$ and $B\in\CD$  that is two points connected by the condition $A=B+m\Go_1+n\Go_2$, $m,n=0\pm 1,\pm 2,\ldots$.
The main vector along any curve  $\GG=A  B$  is defined  by \cite{mus}
\beq
 X+iY=-i[g(z)]_A^B,
\label{3.12}
\eeq
where $[g(z)]_A^B$ is the increment of the function $g(z)$ when $z$ traverses the contour $\GG$ and 
\beq
g(z)=\Gvf(z)+z\ov{\Gvf'(z)}+\ov{\psi(z)}.
\label{3.13}
\eeq
This expression is independent of the contour profile and
equals $g(B)-g(A)$.
The functions $\Gvf(z)$ and $\psi(z)$ are the Muskhelishvili displacement potentials, $\GF(z)=\Gvf'(z)$, $\Psi(z)=\psi'(z)$, and obtained by integrating the series (\ref{3.4})
$$
\Gvf(z)=\Ga_0z-\Ga_1\Gl^2\Gz(z)+\sum_{k=2}^\infty\Ga_{k}\fr{\Gl^{2k}\wp^{(2k-3)}(z)}{(2k-1)!},
$$
\beq
\psi(z)=\Gb_0z-\Gb_1\Gl^2\Gz(z)-\Ga_1\Gl^2\CN(z)+
\sum_{k=2}^\infty\Gb_{k}\fr{\Gl^{2k}\wp^{(2k-3)}(z)}{(2k-1)!}
-\sum_{k=2}^\infty\Ga_{k}\fr{\Gl^{2k}\CN^{(2k-2)}(z)}{(2k-1)!}.
\label{3.14}
\eeq
Since $A$ and $B$ are congruent points, the main vector has to vanish and therefore  $g(B)=g(A)$. In Ref \cite{fil} this condition is considered in the case when
the potentials $\GF(z)$ and $\Psi(z)$ are not only even but also meet the symmetry conditions
\beq
\GF(\bar z)=\ov{\GF(z)}, \quad \GY(\bar z)=\ov{\GY(z)}, \quad z\in\CD,
\label{3.15}
\eeq 
which imply that all the coefficients $\Ga_k$ and $\Gb_k$ are real. In our case, when the angle $\Ga$ is arbitrary or, equivalently, for the chiral vector is not necessarily zigzag or armchair, in general, the conditions (\ref{3.15}) are not valid, and
the coefficients  $\Ga_k$ and $\Gb_k$ are complex. Since the Weierstrass $\Gz$-function satisfies the quasi-periodicity condition 
\beq
\Gz(z+\Go_j)-\Gz(z)=\Gd_j, \quad \Gd_j=2\Gz\left(\fr{\Go_j}{2}\right), \quad j=1,2,
\label{3.16}
\eeq
the function $\Gvf(\Gz)$ is also quasi-doubly-periodic and
\beq
\Gvf(\Gz+\Go_j)-\Gvf(z)=\Ga_0\Go_j-\Ga_1\Gl^2\Gd_j, \quad j=1,2.
\label{3.17}
\eeq
Denote by $\Gc(\Gz)=z\ov{\Gvf'(z)}+\ov{\psi(z)}$.   The relations (\ref{3.8}), (\ref{3.9}) and (\ref{3.14}) imply
\beq
\Gc(z+\Go_j)-\Gc(z)= \bar\Ga_0\Go_j+\bar \Gb_0\bar\Go_j-\bar\Ga_1\Gl^2\bar\Gg_j-\bar\Gb_1\Gl^2\bar\Gd_j, \quad j=1,2.
\label{3.18}
\eeq
On satisfying the condition $g(B)-g(A)=0$ we deduce the following system of two equations:
\beq
2\R\Ga_0\Go_j+\bar\Gb_0\bar\Go_j-\Gl^2(\Gd_j\Ga_1+\bar\Gg_j\bar\Ga_1)-\bar\Gb_1\Gl^2\bar\Gd_j=0, \quad j=1,2.
\label{3.19}
\eeq
Its solution is simplified if the Legendre relation \cite{han} for the case $\R(i\Go_2/\Go_1)<0$ 
\beq
\Gd_1\Go_2-\Gd_2\Go_1=2\pi i
\label{3.20}
\eeq
is employed. We have
$$
\R\Ga_0=\fr{1}{2 \I (\Go_1\bar\Go_2)}[\Gl^2\I\{\Ga_1(\Gd_1\bar\Go_2-\Gd_2\bar\Go_1)\}, 
$$
\beq
\bar\Gb_0=\fr{1}{2 \I (\Go_1\bar\Go_2)}\left[\pi\Gl^2\Ga_1-\fr{i}{2}\Gl^2\bar\Ga_1(\bar\Gg_1\Go_2-\bar\Gg_2\Go_1)-\fr{i}{2}\Gl^2\bar\Gb_1(\bar\Gd_1\Go_2-\bar\Gd_2\Go_1)\right].
\label{3.21}
\eeq
where for the hexagonal cell under consideration $ \I (\Go_1\bar\Go_2)=-\fr{a^2\sqrt{3}}{2}$.
Formulas (\ref{3.21}) can be further simplified if we use the identities 
\beq
\Gd_1\bar\Go_2-\Gd_2\bar\Go_1=0
\label{3.22}
\eeq
valid
for  hexagonal lattices. To prove them, note that for this type of  symmetry
\beq
\Gz(e^{i\pi/6})=e^{-i\pi/6}\left\{\fr{1}{z}+\sum_{m,n}^{} {}'\left[\fr{1}{z-we^{-i\pi/6}}+\fr{1}{we^{-i\pi/6}}+\fr{z}{(w e^{-i\pi/6})^2}\right]\right\}.
\label{3.23}
\eeq
Therefore 
\beq
\Gz(e^{i\pi/6})=e^{-i\pi/6}\Gz(z), \quad \Gz\left(\fr{\Go_2}{2}\right)= e^{-i\pi/3}\Gz\left(\fr{\Go_1}{2}\right), 
\label{3.24}
\eeq
and the identity (\ref{3.22}) follows.
Notice the relation (\ref{3.22}) is known \cite{fil}  for equilateral triangular and square lattices. The identity (\ref{3.22}) immediately implies that the constants in (\ref{3.9'}) vanish,
$\Gg_1=\Gg_2=0$. Using these results we find the final formulas for the coefficients $\R\Ga_0$ and $\Gb_0$
\beq
\R\Ga_0=-\fr{\pi\Gl^2\Gb_1}{a^2\sqrt{3}}, \quad
\Gb_0=-\fr{2\pi\Gl^2\bar\Ga_1}{a^2\sqrt{3}}.
\label{3.25}
\eeq
It follows from here that $\I\Gb_1=0$. It will shown later that the solution is independent of $\I\Ga_0$.

\subsection{Infinite system of linear algebraic equations}

The series representations (\ref{3.4}) of the functions $\GF(z)$ and $\Psi(z)$ have arbitrary complex coefficients
$\Ga_k$ and $\Gb_k$. They have to be fixed from the boundary condition of the circle $L_{00}$ 
(\ref{2.6}) that can be written in the form \cite{eng}
\beq
\GF(t)+\ov{\GF(t)}-[\bar t\GF'(t)+\Psi(t)]e^{2i\Gt}=-\Gs_+-\Gs_-e^{2i(\Gt-\Ga)}, \quad 0\le\Gt\le 2\pi.
\label{3.26}
\eeq
In what follows we aim to recast the boundary condition as an infinite system for the unknown coefficients.
 
 By representing the Weierstrass function $\wp(z)$ by its Laurent series \cite{han}
\beq
\wp(z)=\fr{1}{z^2}  +c_2z^2+\ldots+c_s z^{2s-2}+\ldots, \quad z\in\CA,
\label{3.27}
\eeq
convergent uniformly and absolutely in any annulus $\CA=\{0<r_-\le|z|\le r_+<\fr{a}{2}\}$. The coefficients $c_s$ are expressed through  the double series
\beq
c_s=(2s-1)\sum_{m,n}^{} {}'\fr{1}{w^{2s}}, \quad  w=m\Go_1+n\Go_2, \quad m,n=0,\pm 1,\pm 2,\ldots.
\label{3.28}
\eeq
In addition, the coefficients $c_s$ can be recovered by the recursion formula
\beq
c_s=\fr{3}{(2s+1)(s-3)}\sum_{t=2}^{s-2}c_t c_{s-t}, \quad s>3,
\label{3.29}
\eeq
after the coefficients $c_2$  and $c_3$ are computed by (\ref{3.28}).  By differentiating $2k$-times the Laurent series (\ref{3.27}),
\beq
\wp^{(2k)}(z)=\fr{(2k+1)!}{z^{2k+2}}+\sum_{j=0}^\infty
\fr{(2k+2j)!c_{k+j+1}z^{2j}}{(2j)!},
\quad k=1,2,\ldots.
\label{3.30}
\eeq
The Laurent series of the derivatives of the Natanzon function has a similar form \cite{fil}
$$
\CN'(z)=\sum_{j=1}^\infty (2+2j)(1+2j) d_{j+1} z^{2j}, 
$$
\beq
\CN^{(2k+1)}(z)=\sum_{j=0}^\infty\fr{(2k+2+2j)! d_{k+j+1} z^{2j}}
{(2j)!}, \quad k=1,2,\ldots,
\label{3.31}
\eeq
with the coefficients $d_s$ defined by 
\beq
d_s=\sum_{m,n}^{} {}'\fr{\bar w}{w^{2s+1}}, \quad s=2,3,\ldots.
\label{3.32}
\eeq
For these coefficients, there is no analog of the recursion formula (\ref{3.29}) for  the function  $\CN(z)$ is not doubly-periodic and 
has no analog of the differential equation 
\beq
[\wp'(z)]^2-4\wp^3(z)+g_2\wp(z)+g_3=0, \quad 
g_2=60\sum_{m,n}^{} {}'\fr{1}{w^4}, \quad g_3=140\sum_{m,n}^{} {}'\fr{1}{w^6},
\label{3.34}
\eeq 
of the Weierstrass function that is exploited for the derivation of the recursion formula (\ref{3.29}) for the coefficients $c_s$.
The expressions for the coefficients $c_s$ and $d_s$ can be simplified if the symmetry of the hexagonal lattice is taking into account.
Similar to the regular triangular lattice \cite{fil} for the case of hexagonal symmetry we have
\beq
c_{3s-1}=c_{3s-2}=0, \quad d_{3s}=d_{3s+1}=0, \quad s=1,2,\ldots.
\label{3.35}
\eeq
The nonzero coefficients are computed by
\beq
c_{3s}=(6s-1)\sum_{m,n}^{} {}'\fr{1}{w^{6s}}, \quad d_{3s-1}=\sum_{m,n}^{} {}'\fr{\bar w}{w^{6s-1}}, 
 \quad s\ge 1.
\label{3.37}
\eeq
The recursion formula (\ref{3.29}), when adjusted to hexagonal case, has the form 
\beq
 c_{3s}=\fr{1}{(6s+1)(s-1)}\sum_{t=1}^{s-1}c_{3t} c_{3(s-t)}, \quad s\ge 2.
\label{3.38}
\eeq
Substituting  the series expansions of the derivatives of the functions $\wp(z)$ and $\CN(z)$ into formulas (\ref{3.4}) yields \cite{fil}
$$
\GF(z)=\Ga_0+\sum_{k=0}^\infty\Ga_{k+1}\Gl^{2k+2}\left(\fr{1}{z^{2k+2}}+\sum_{j=0}^\infty r_{j,k}z^{2j}\right),
$$
\beq
\GY(z)=\Gb_0+\sum_{k=0}^\infty\Gb_{k+1}\Gl^{2k+2}\left(\fr{1}{z^{2k+2}}+\sum_{j=0}^\infty r_{j,k}z^{2j}\right)
-\sum_{k=0}^\infty\Ga_{k+1}\Gl^{2k+2}\sum_{j=0}^\infty\Gr_{jk}z^{2j}, 
\label{3.39}
\eeq
where
\beq
r_{jk}=\fr{(2k+2j)! c_{j+k+1}}{(2k+1)!(2j)!}.
\quad \Gr_{jk}=\fr{(2k+2+2j)! d_{j+k+1}}{(2k+1)!(2j)!}, \quad j,k=0,1,\ldots,
\label{3.40}
\eeq
and $r_{00}=\Gr_{00}=0$.
On the boundary $L_{00}$, $z=\Gl e^{i\Gt}$. 
Changing the order of summation in (\ref{3.39} and comparing the same powers of $e^{i\Gt}$
in the boundary condition (\ref{3.26}) enable us to deduce the following system for the unknown coefficients:
$$
2\R\Ga_0-\Gb_1+2\sum_{k=0}^\infty \Gl^{2k+2} r_{0k}\R \Ga_{k+1}=-\Gs_+,
$$
$$
\Gb_{j+1}=(2j+1)\Ga_j+\sum_{k=0}^\infty \Gl^{2k+2j+2} r_{jk}\ov{\Ga_{k+1}}, \quad j=1,2,\dots,
$$
\beq
\ov{\Ga_j}+\sum_{k=0}^\infty \Gl^{2k+2j}\{[\Gl^2(1-2j)r_{jk}+\Gr_{j-1k}]\Ga_{k+1}-r_{j-1 k}\Gb_{k+1}\}=(\Gb_0-\Gs_-e^{-2\Ga i})\Gd_{j1}, \quad j=1,2,\dots.
\label{3.41}
\eeq
We next represent the complex coefficients as $\Ga_k=\Ga_k'+i\Ga_k''$, express $\Gb_k$ from the second equation in (\ref{3.41}), use the relations 
(\ref{3.25}) for $\R\Ga_0$ and $\Gb_0$ and obtain the following two separate systems of real equations for the real and imaginary  parts of the coefficients $\Ga_k$:
$$
(b-1)\Gb_1+2\sum_{k=0}^\infty \Gl^{2k+2} r_{0k}\Ga'_{k+1}=-\Gs_+,
$$
\beq
\Ga_j'+\sum_{k=1}^\infty\Gl^{2k+2j}\left(d^-_{jk}+\fr{2}{b-1}r_{j-1, 0}r_{0, k-1}\right)\Ga_k'-b\Gd_{j1}\Ga_1'=-\fr{\Gs_+\Gl^{2j}}{b-1}r_{j-1, 0}-\Gs_-\Gd_{j 1}\cos 2\Ga,
\label{3.42}
\eeq
and 
\beq
\Ga_j''-\sum_{k=1}^\infty\Gl^{2k+2j}d^+_{jk}\Ga_k''+b\Gd_{j1}\Ga_1''=-\Gs_-\Gd_{j 1}\sin 2\Ga, \quad j=1,2,\ldots.
\label{3.43}
\eeq
Here,
$$
b=\fr{2\pi\Gl^2}{\sqrt{3} a^2},
$$
\beq
d^\pm_{jk}=(1-2j)r_{j,k-1}-(1+2k)r_{j-1,k}+\Gl^{-2}\Gr_{j-1.k-1}
\pm\sum_{m=1}^\infty \Gl^{4m}r_{j-1,m}r_{m,k-1}.
\label{3.44}
\eeq

After the systems (\ref{3.42}) and  (\ref{3.43}) for the coefficients $\Ga_k=\Ga_k'+i\Ga_k''$ are solved, the coefficient $\Gb_1$ is determined by the first equation in 
(\ref{3.42}), the other coefficients  $\Gb_k$ $(k=2,3,\ldots)$ are expressed through $\Ga_k$
by the second equation in (\ref{3.41}), and the coefficients $\Ga_0'=\Re \Ga_0$ and $\Gb_0$ are fixed by equations (\ref{3.25}).

\setcounter{equation}{0}

\section{Stresses, displacements, and elastic constants of the graphene bonds}

\begin{figure}[t]
\centerline{
\scalebox{0.6}{\includegraphics{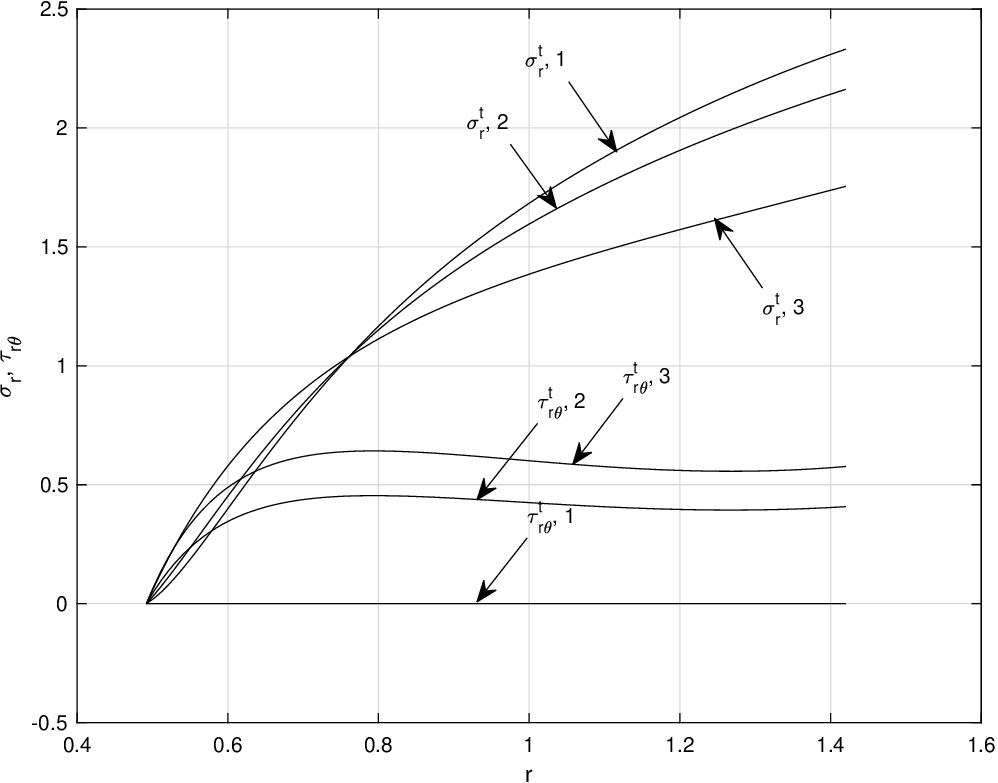}}}
\caption{
The total stresses $\Gs^t_r$ and $\tau^t_{r\theta}$ vs $r\in[\Gl, \fr{a}{\sqrt{3}}]$ for 
 $\Ga=0$ (curves 1),  $\Ga=\fr{\pi}{8}$ (curves 2),  and $\Ga=\fr{\pi}{4}$ (curves 3) when
$\Gl=\fr{a}{5}$,  $\Gs_1=2$, $\Gs_2=1$,   $\Gt=0$.}
\label{fig2}
\end{figure} 

\begin{figure}[t]
\centerline{
\scalebox{0.6}{\includegraphics{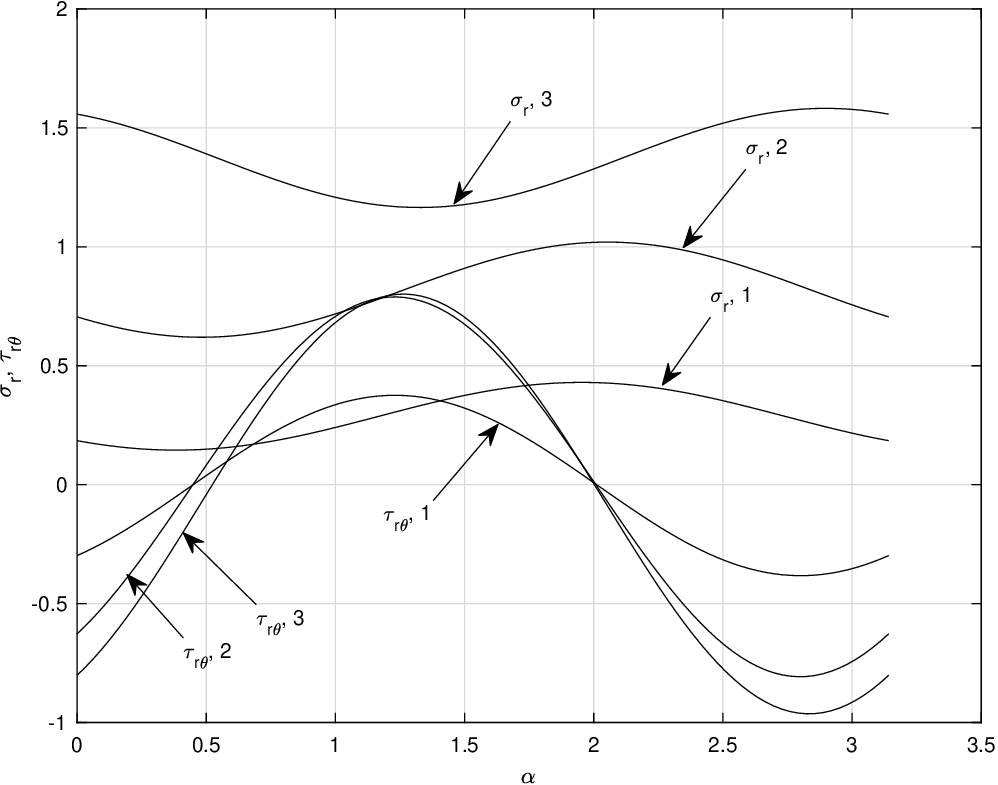}}}
\caption{
The total stresses $\Gs^t_r$ and $\tau^t_{r\theta}$ vs $\Ga\in[0,\pi]$ for $\Gt=\fr{\pi}{8}$ and
$r=\Gl+\fr{1}{16}(\fr{a}{2}-\Gl)$ (curves 1), $r=\Gl+\fr{1}{4}(\fr{a}{5}-\Gl)$ (curves 2), and $r=\Gl+\fr{15}{16}(\fr{a}{2}-\Gl)$
(curves 3) when
$\Gl=\fr{a}{5}$,  $\Gs_1=2$, and $\Gs_2=1$.}
\label{fig3}
\end{figure}

The total stresses $\Gs^t_r$
and $\tau^t_{r\Gt}$ are recovered by
 \beq
 \Gs^t_r-i\tau^t_{r\Gt}=\Gs_r^{K}-i\tau_{r\Gt}^{K}+2\R\GF(\tau)-[\bar\tau\GF'(\tau)+\Psi(\tau)]e^{2i\Gt}, \quad \tau=re^{i\Gt}, \quad \Gl\le r\le\fr{a}{\sin\Gt+\sqrt{3}\cos\Gt},
 \label{4.1}
 \eeq
 where $\Gs_r^K$ and $\tau_{r\Gt}^{K}$ are given by (\ref{2.2}). The functions $\GF(\tau)$ and $\Psi(\tau)$ are determined by series (\ref{3.39}), while the derivative $\GF'(\tau)$
 is computed by
 \beq
\GF'(\tau)=\sum_{k=0}^\infty\Ga_{k+1}\Gl^{2k+2}\left(-\fr{2k+2}{\tau^{2k+3}}+2\sum_{j=1}^\infty jr_{j,k}\tau^{2j-1}\right).
\label{4.1'}
\eeq
 Numerical results  for the normal and tangential stresses $\Gs_r$ and $\tau_{r\Gt}$ in polar coordinates $(r,\Gt)$, $\Gt=0$, versus $r$, 
 $\Gl\le r\le \fr{a}{\sqrt{3}}$,
 are shown in Fig. 2. The remote loads are $\Gs_1=2$ and $\Gs_2=1$, and the load $\Gs_1$ forms angle $\Ga$ with the $x$-axis.
 For curves $1$, $2$, and $3$, the angle $\Ga$ is chosen to be $0$, $\fr{\pi}{8},$ and $\fr{\pi}{4}$, respectively. 
It is experimentally confirmed that an atom of hydrogen  cannot pass through dense clouds of electrons
forming the bonds of a graphene cell \footnote{K. S. Novoselov, private communication}.  The radius of an atom of hydrogen is $\Gr=53$ pm,
while for the graphene cell $\half a=123$ pm, and the ratio $\fr{\Gr}{a/2}=0.431$. Based on this argument, the hole radius $\Gl\le \Gr$,
and for computations we select $\Gl=0.4 \fr{a}{2}=49$ pm. It is seen from Fig.2 that for $\Gt=0$ the normal stress $\Gs_r$ attains
its maximum for $r=\fr{a}{\sqrt{3}}$. All the curves are initiated at the point $(\Gl,0)$ that is the stresses vanish at 
the hole boundary, and the boundary conditions are fulfilled.

 Variation of the normal and tangential stresses $\Gs_r$ and $\tau_{r\Gt}$ with the angle $\Ga$ is illustrated 
 in Fig. 3. The polar angle $\Gt=\fr{\pi}{8}$, while $r=\Gl+\fr{1}{16}(\fr{a}{2}-\Gl)$ (curves 1), $r=\Gl+\fr{1}{4}(\fr{a}{5}-\Gl)$ (curves 2), and $r=\Gl+\fr{15}{16}(\fr{a}{2}-\Gl)$. It is seen that the amplitude of the stresses increases as $r$ grows and varies significantly with
 variation of the angle $\Ga$ or, equivalently, the chiral vector $C_h=(m,n)$.

 The  total displacements $u^t$ and $v^t$ are expressed through the remote loads and the displacement potentials 
 \beq
 2G(u^t+iv^t)=\fr{\Gk-1}{4}(\Gs_1+\Gs_2)z+\fr{1}{2}(\Gs_1-\Gs_2)\bar z+h(z),
 \label{4.2}
 \eeq
 where
 \beq
 h(z)=\Gk\Gvf(z)-z\ov{\Gvf'(z)}-\ov{\psi(z)}.
 \label{4.3}
 \eeq
 The parameter $\Gk$  associated with the generalized plane strain state is
 \beq
 \Gk=\fr{3-\nu}{1+\nu},
 \label{4.4}
 \eeq
 and  $G$ and $\nu$ are  the shear modulus and Poisson ratio of the graphene bonds.
The functions $(\ref{4.3})$ are expanded in terms of odd derivatives of the Weierstrass function $\wp(z)$ and even derivatives of the Nantanzon function $\CN(z)$ in (\ref{3.14}).
On differentiating the Laurent series (\ref{3.30}) and (\ref{3.31}) we eventually find
the series expansion of the displacement $(u^t+iv^t)$
$$
2G(u^t+iv^t)=\fr{\Gk-1}{4}(\Gs_1+\Gs_2)z+\fr{1}{2}(\Gs_1-\Gs_2)\bar z+\Ga'_0 (\Gk-1) z-\bar\Gb_0\bar z
$$
$$
+\Gk\sum_{k=1}^\infty\Ga_{2k}\Gl^{2k}\left(\fr{z^{1-2k}}{1-2k}+\sum_{j=0}^\infty r_{j,k-1}\fr{z^{1+2j}}{1+2j}\right)-z\sum_{k=1}^\infty\bar\Ga_{2k}\Gl^{2k}
\left(\fr{1}{\bar z^{2k}}+\sum_{j=0}^\infty r_{j,k-1}\bar z^{2j}\right)
$$
\beq
-\sum_{k=1}^\infty\bar\Gb_{2k}\Gl^{2k}\left(\fr{\bar z^{1-2k}}{1-2k}+\sum_{j=0}^\infty r_{j,k-1}\fr{\bar z^{1+2j}}{1+2j}\right)+\sum_{k=1}^\infty\bar\Ga_{2k}\Gl^{2k}
\sum_{j=0}^\infty \Gr_{j,k-1} \fr{\bar z^{1+2j}}{1+2j}.
\label{4.5}
\eeq

This formula  has two unknowns, the shear modulus $G=\fr12(1+\nu)^{-1}E$ and the parameter $\Gk$ expressed by (\ref{4.4}) through
the Poisson ratio $\nu$.  The Young modulus $E$ and Poisson ratio $\nu$ of the graphene bonds
do not coincide with their effective counterparts $E^\circ$ and $\nu^\circ$ measured by laboratory tests.
To determine $E$ and  $\nu$ of the graphene bonds, we apply the reverse version of the method proposed \cite{fil} for recovering the effective constants
for a general doubly-periodic structures and specified for triangular and square lattices.
The main idea of the method is that the displacement in the homogenized structure and the perforated plane have to be quasi-doubly-periodic functions
with the same lines of discontinuities and the same jumps. As in the case of triangular and square lattices, the homogenized plane associated with the hexagonal doubly-periodic 
perforated plane is isotropic. To show this assume that it is orthotropic with $ E^\circ_1$ and $\nu^\circ_1$ being the Young modulus and Poisson ratio in the $x$-direction (armchair)
and
$E^\circ_2$ and $\nu^\circ_2$ being the corresponding constants in the $y$-direction (zigzag).  By writing the Hooke's law equations for the homogenized plane
subjected to the remote loads $\Gs_1$ and $\Gs_2$ with $\Ga=0$ as in Fig. 1 and using  the fact that the stresses everywhere in the plane are constants, $\Gs_x=\Gs_1$, $\Gs_y=\Gs_2$, $\tau_{xy}=0$, we have
$$
\fr{\Md u^\circ}{\Md x}=\fr{\Gs_1}{E^\circ_1}-\fr{\nu^\circ_2\Gs_2}{ E^\circ_2},\quad \fr{\Md u^\circ}{\Md y}+\fr{\Md v^\circ}{\Md x}=0, 
$$
\beq
\fr{\Md v^\circ}{\Md y}=\fr{\Gs_2}{ E^\circ_2}-\fr{\nu^\circ_1\Gs_1}{E^\circ_1},\quad \nu^\circ_2E^\circ_1=\nu^\circ_1E^\circ_2.
\label{4.6}
\eeq
We integrate these  relations and obtain the displacements. Up to the rigid body displacements $(c_x,c_y)$ and rotation $\Go$, they are
\beq
u^\circ+iv^\circ=(\Gk_1^-\Gs_1+\Gk_2^-\Gs_2)z+(\Gk_1^+\Gs_1-\Gk_2^+\Gs_2)\bar z +c_x+i c_y-i\Go z,
\label{4.7}
\eeq
where 
\beq
\Gk_j^\pm=\fr{1\pm \nu^\circ_j}{E^\circ_j}, \quad j=1,2.
\label{4.8}
\eeq
From here, the doubly-quasi-periodicity relations for the effective displacements immediately follow. By disregarding the rigid body rotation we obtain
$$
(u^\circ+iv^\circ)(z+\Go_j)-(u^\circ+iv^\circ)(z)
$$
\beq
=[(\Gk_1^-+\Gk_2^-)\Go_j+(\Gk_1^+-\Gk_2^+)\bar\Go_j]\Gs_++
[(\Gk_1^--\Gk_2^-)\Go_j+(\Gk_1^++\Gk_2^+)\bar\Go_j]\Gs_-, \quad j=1,2.
\label{4.9}
\eeq
On the other side, if we use formulas (\ref{4.2}), (\ref{4.3}) and also (\ref{3.17}) and (\ref{3.18}) we derive
the analog of the relations (\ref{4.9}) for the total displacements of the graphene sheet. All  coefficients $\Ga_j$ and $\Gb_j$ are real now, and we have
$$
(u^t+iv^t)(z+\Go_j)-(u^t+iv^t)(z)
$$
\beq
=\fr{1-\nu}{E}\Gs_+\Go_j+\fr{1+\nu}{E}\Gs_-\bar\Go_j
+\fr{1+\nu}{E}[\Ga_0(\Gk-1)\Go_j-\Gb_0\bar\Go_j-\Ga_1\Gl^2\Gk\Gd_j+\Gb_1\Gl^2\bar\Gd_j], \quad j=1,2.
\label{4.10}
\eeq
We can represent the displacement jumps of the homogenized plane and the perforated plane in the basis $\{\Gs_+,\Gs_-\}$ and
denote $\Ga_j=\Ga_j^+$, $\Gb_j=\Gb_j^+$ 
for $\Gs_+=1$, $\Gs_-=0$ and $\Ga_j=\Ga_j^-$, $\Gb_j=\Gb_j^-$ 
for $\Gs_+=0$, $\Gs_-=1$.  It is directly verified that
\beq
\Ga_1^+=\Gb_0^+=0,\quad \Ga_0^-=\Gb_1^-=0.
\label{4.11}
\eeq 
On equating the right hand-sides of equations (\ref{4.9}) and (\ref{4.10}) for these two sets of solutions we arrive at the 
system of four equations 
$$
\Gk_1^-\Go_1+\Gk_2^-\Go_1+\Gk_1^+\bar\Go_1-\Gk_2^+\bar\Go_1=b^+_1,
$$
$$
\Gk_1^-\Go_2+\Gk_2^-\Go_2+\Gk_1^+\bar\Go_2-\Gk_2^+\bar\Go_2=b^+_2,
$$
$$
-\Gk_1^-\Go_1+\Gk_2^-\Go_1-\Gk_1^+\bar\Go_1-\Gk_2^+\bar\Go_1=b^-_1,
$$
\beq
-\Gk_1^-\Go_2+\Gk_2^-\Go_2-\Gk_1^+\bar\Go_2-\Gk_2^+\bar\Go_2=b^-_2,
\label{4.12}
\eeq
where
$$
b_j^+=
\fr{1-\nu}{E}\Go_j+\fr{1+\nu}{E}[\Ga_0^+(\Gk-1)\Go_j+\Gb_1^+\Gl^2\bar\Gd_j],
$$
\beq
b_j^-=
\fr{1+\nu}{E}(-\bar\Go_j+\Ga_1^-\Gl^2\Gk\Gd_j+\Gb_0^-\bar\Go_j),\quad j=1,2.
\label{4.13}
\eeq
The determinant $\GD$ of the system with respect to $\Gk_j^\pm$ is 
\beq
\GD=4(\Go_2\bar\Go_1-\Go_1\bar\Go_2)^2=-12 a^4,
\label{4.14}
\eeq
and the system has a unique solution. 
Notice that since $\Go_1=\bar\Go_2$, the complex numbers $\Gd_j$ introduced in (\ref{3.11}) in the case of a hexagonal lattice,  satisfy
the relations
\beq
\Gd_1\bar\Go_2-\Gd_2\bar\Go_1=0, \quad \Gd_1=\bar\Gd_2, \quad \fr{\bar\Gd_1}{\Go_1}=\fr{\bar\Gd_2}{\Go_2}=\Gd,
\label{4.15}
\eeq
where
$\Gd$ is a real number. That is why the solution of the system (\ref{4.12}) meets the conditions
\beq
\Gk_1^+=\Gk_2^+=\Gk^+=\fr{1+\nu^\circ}{E^\circ}, \quad \Gk_1^-=\Gk_2^-=\Gk^-=\fr{1-\nu^\circ}{E^\circ}.
\label{4.16}
\eeq
This implies that the homogenized graphene is an isotropic material with the Poisson ratio $\nu^\circ$ and the Young modulus $ E^\circ$.
We recall that $\Gk$ is given by (\ref{4.4}) and rewrite the system (\ref{4.12}) as a system of two equations with respect to $1/E$ and $\nu/E$
$$
(1+2\Ga_0^++\Gb_1^+\Gl^2\Gd)\fr{1}{E}+(-1-2\Ga_0^++\Gb_1^+\Gl^2\Gd)\fr{\nu}{E}=\fr{1-\nu^\circ}{E^\circ},
$$
\beq
(1-\Gb_0^- -3\Ga_1^-\Gl^2\Gd)\fr{1}{E}+(1-\Gb_0^-+\Ga_1^-\Gl^2\Gd)\fr{\nu}{E}=\fr{1+\nu^\circ}{E^\circ}.
\label{4.17}
\eeq
\begin{figure}[t]
\centerline{
\scalebox{0.6}{\includegraphics{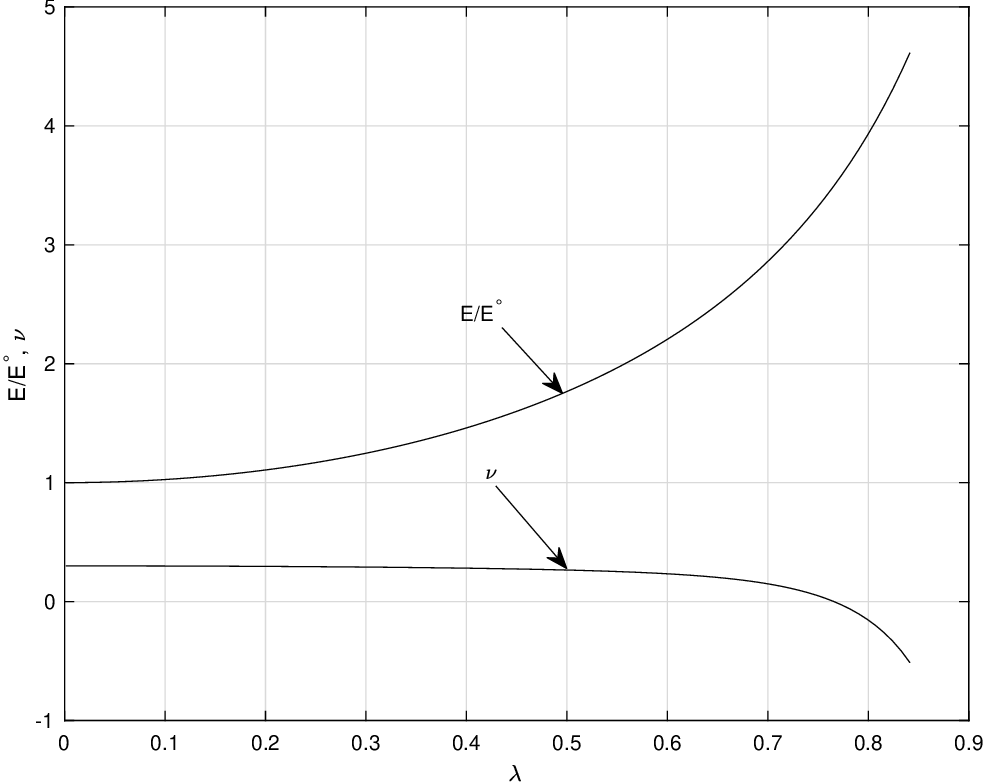}}}
\caption{Poisson ratio $\nu$ and the dimensionless parameter $E/E^\circ$ of the bonds vs radius $\Gl$ when  $\nu^\circ=0.3$ and $ E^\circ$ is fixed.}
\label{fig4}
\end{figure}

On solving it we find the expressions of the Poisson ratio and the Young modulus through their effective counterparts
\beq
\nu=\fr{\GD_1}{\GD_0},\quad 
\fr{E}{E^\circ}=\fr{2\GD_2}{\GD_0},
\label{4.18}
\eeq
where
$$
\GD_0=2-(1-\nu^\circ)(\Gb_0^--\Ga_1^-\Gl^2\Gd)-(1+\nu^\circ)(\Gb_1^+\Gl^2\Gd-2\Ga_0^+)
$$
$$
\GD_1=2\nu^\circ+(1-\nu^\circ)(\Gb_0^-+3\Ga_1^-\Gl^2\Gd)+(1+\nu^\circ)(\Gb_1^+\Gl^2\Gd+2\Ga_0^+)
$$
\beq
\GD_2=2\Ga_1^-\Gb_1^+\Gl^4\Gd^2-(1+2\Ga_0^+)(\Ga_1^-\Gl^2\Gd+\Gb_0^- -1).
\label{4.19}
\eeq

In Fig. 4 we plot the dimensionless parameter $E/E^\circ$  and the Poisson ratio $\nu$ as functions of the hole radius $\Gl$
when the effective Poisson ratio is fixed, $\nu^\circ=0.3$. When bonds thickness decreases the Young modulus is increasing,
while Poisson ratio decreases. For the value $\Gl=\fr{a}{5}$ selected for Figs 2 and 3, the values of the bond constants are
$\nu=0.2668$ and $E=1.7377E^\circ$ that is the Young modulus exceeds the value of its effective counterpart by 74\%.

\begin{figure}[t]
\centerline{
\scalebox{0.6}{\includegraphics{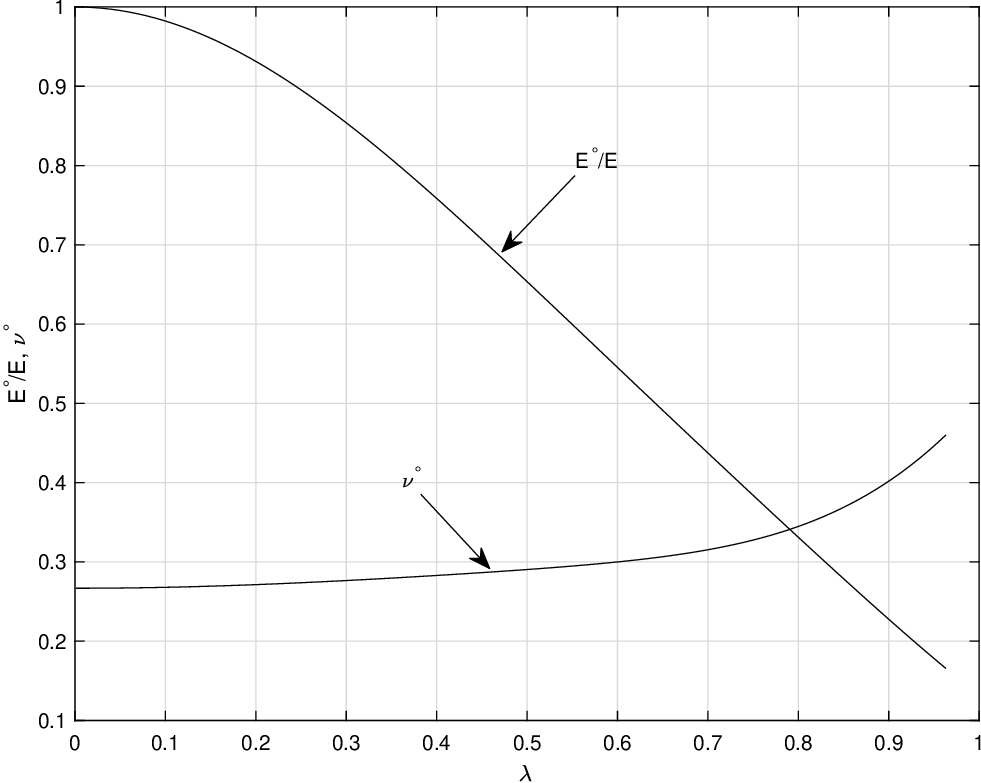}}}
\caption{Effective Poisson ratio $\nu^\circ$ and the dimensionless parameter $E^\circ/ E$ of the perforated plane vs radius $\Gl$ when  $\nu=0.2668$ and $E$ is fixed.}
\label{fig5}
\end{figure} 

It is of interest to express the effective constants $E^\circ$ and $\nu^\circ$ the bond constants $E$ 
and $\nu$. By solving the system (\ref{4.17}) with respect to $E^\circ$ and $\nu^\circ$ we discover
\beq
\fr{E^\circ}{E}=\fr{2}{1-\nu+(1+\nu)(\GL^-+\GL^+)}, \quad 
\nu^\circ=\fr{\nu-1+(1-\nu)(\GL^- -\GL^+)}{1-\nu+(1+\nu)(\GL^-+\GL^+)},
\label{4.20}
\eeq
where 
\beq
\GL^+=\Ga_0^+(\Gk-1)+\Gb_1^+\Gl^2\Gd, \quad
\GL^-=1-\Ga_1^-\Gl^2\Gd\Gk-\Gb_0^-,
\label{4.21}
\eeq
and $\Gk$ is the parameter defined in (\ref{4.4}). 
The results of computations are shown for the Poisson ratio $\nu=0.2688$.
The effective Young modulus $E^\circ$ decreases with the growth of the 
hole radius $\Gl$, while the effective Poisson ratio $\nu^\circ$ increases.

\begin{figure}[t]
\centerline{
\scalebox{0.6}{\includegraphics{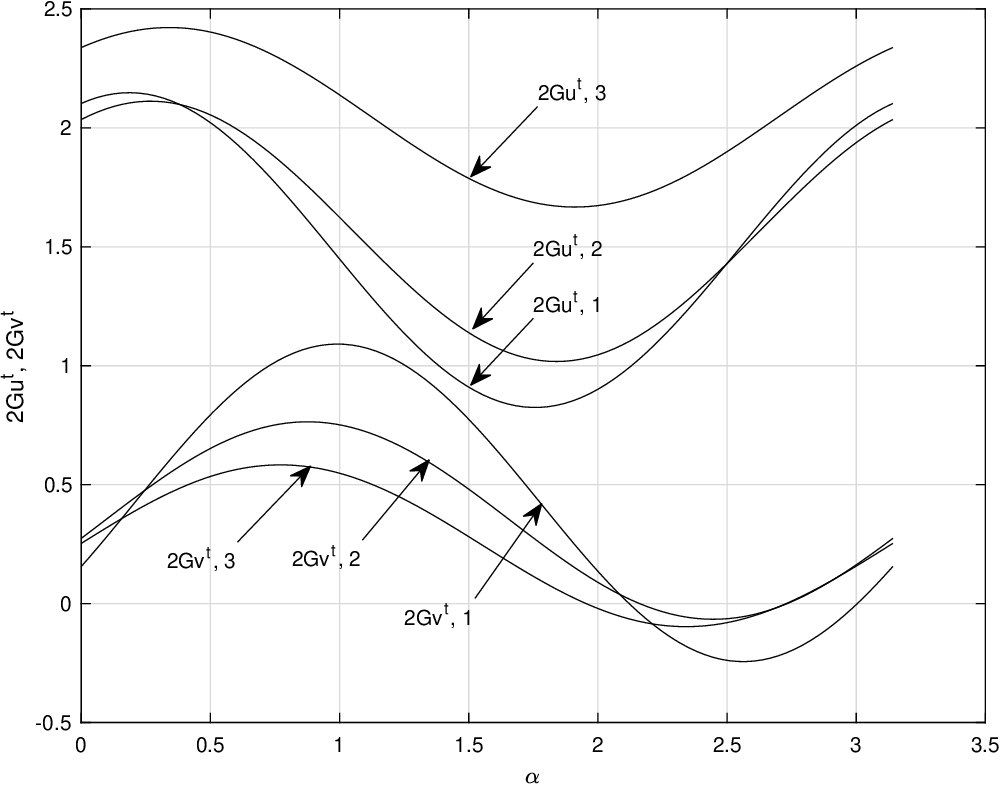}}}
\caption{
The total displacements $2Gu^t$ and $2Gv^t$ vs $\Ga\in[0,\pi]$ for $\Gt=\fr{\pi}{8}$ and 
$r=\Gl$ (curves 1), $r=\Gl+\fr{1}{3}(\fr{a}{2}-\Gl)$ (curves 2), and $r=\fr{a}{2}$ (curves 3) when
$\Gl=\fr{a}{5}$,  $\Gs_1=2$, $\Gs_2=1$, and $\nu=0.2668$.}
\label{fig6}
\end{figure} 

For computation of the displacements by formula (\ref{4.5}) we use the  value $\Gl=\fr{a}{5}$ for the hole radius and the corresponding
value of the Poisson ratio $\nu=0.2688$ computed by formula (\ref{4.18}). In Fig. 6 we show the total displacements $2Gu^t$ and $2G v^t$ 
as functions of the angle $\Ga\in [0,\pi]$ for $r=\Gl$ (curves 1), $r=\Gl+\fr{1}{3}(\fr{a}{2}-\Gl)$ (curves 2), and $r=\fr{a}{2}$
(curves 3). The other parameters are $\Gt=\fr{\pi}{8}$, $\Gs_1=2$, and $\Gs_2=1$. It is seen that the highest variation of the displacements 
occurs on the hole boundary.

\vspace{.2in}

{\large {\bf Conclusions}}

\vspace{.1in}

The paper presents a plane stress model of graphene subjected to tension at infinity. A graphene sheet characterized by a general chiral vector
is modeled as a
doubly-periodic thin plate whose unit cell is a hexagon with a central circular hole.  For the solution, the method  of doubly-periodic 
Kolosov-Muskhelishvili complex potentials, the Weierstrass elliptic and Natanzon quasi-doubly-periodic meromorphic functions proposed for the 
armchair-zigzag symmetry \cite{fil} is generalized to any chiral vector case. In contrast to the symmetric case, in the general case, the two Kolosov-Muskhelishvili potentials
do not satisfy the symmetry condition $\GF(\bar z)=\ov{\GF(z)}$,  $\GY(\bar z)=\ov{\GY(z)}$, the unit cell is not reducible to an equilateral triangular cell,
and the resulting infinite system is not real and has complex coefficients. 
 
 We have derived series representations of the stresses and displacements. Numerical calculations show that for all chiral vectors tested that
 is for all $\Ga\in[0,\fr{\pi}{6}]$,  
 the radial stress $\Gs_r$ attains it maximum at the vertex $x=\fr{a}{\sqrt{3}}$, $y=0$. Since the displacements comprise a part of the solution to  a plane elasticity
 problem, they depends on the Young modulus $E$ and Poisson ratio $\nu$ of the material or, equivalently in the framework of the continuum mechanics model of graphene,
 the graphene bonds. These elastic moduli are not available in the literature. What are  known by experimental measurements and {\it ab initio} calculations are the effective Young modulus
 $E^\circ\sim 1 $TPa and the effective Poisson ratio $\nu^\circ\in[0.12, 0.413]$. We have derived expressions of the elastic constants of the graphene bonds as 
 functions of their effective counterparts. In particularly, it has been found that for the effective Poisson ratio $\nu^\circ=0.3$ and the radius  $\Gl=\fr{a}{5}$
 of the hole in a unit cell, the Young modulus of the graphene bonds exceeds its effective value $1 TPa$ by $74\%$.

\end{document}